\documentclass{article}
\usepackage{spconf,amsmath,graphicx}
\usepackage{amssymb}
\usepackage{amsmath}
\usepackage{graphicx}
\usepackage{graphicx}
\usepackage{makecell}
\usepackage{float}
\newcommand{\tensor}[1]{\boldsymbol{\mathcal{#1}}}

\newcommand{\mat}[1]{\mathbf{#1}}
\newcommand{\vect}[1]{\mathbf{#1}}

\title{High-order Tensor Completion for Data Recovery via Sparse Tensor-train Optimization}
\name{Longhao Yuan$^{1, 2}$, Qibin Zhao$^{2, 3, *\thanks{*Corresponding authors: qibin.zhao@riken.jp, cao@sit.ac.jp}}$ and Jianting Cao$^{4, 1, *}$}
\address{$^1$Graduate School of Engineering, Saitama Institute of Technology, Japan \and $^2$Tensor Learning Unit, RIKEN Center for Advanced Intelligence Project (AIP), Japan
\and $^3$School of Automation, Guangdong University of Technology, China \and $^4$School of Computer Science and Technology, Hangzhou Dianzi University, China\\}
\begin{document}
\topmargin=0mm
%
\maketitle
\begin{abstract}
In this paper, we aim at the problem of tensor data completion. Tensor-train decomposition is adopted because of its powerful representation ability and linear scalability to tensor order. We propose an algorithm named Sparse Tensor-train Optimization (STTO) which considers incomplete data as sparse tensor and uses first-order optimization method to find the factors of tensor-train decomposition. Our algorithm is shown to perform well in simulation experiments at both low-order cases and high-order cases. We also employ a tensorization method to transform data to a higher-order form to enhance the performance of our algorithm. The results of image recovery experiments in various cases manifest that our method outperforms other completion algorithms. Especially when the missing rate is very high, e.g., 90\% to 99\%, our method is  significantly better than the state-of-the-art methods. 
\end{abstract}
\section{Introduction}
Tensors are multi-dimensional arrays and high-order generation of vectors and matrices \cite{kolda2009tensor}. Most of the real world data like color images, videos, multichannel electroencephalography (EEG) signals, etc. are more than two dimensions. Tensor data representation can keep the original form of data, which is good for retaining high dimensional structure and adjacent relation information of data. Due to the flexibility and highly compressibility of tensor decomposition, in recent decades, many tensor methodologies have been proposed in various fields such as image and video completion \cite{acar2011scalable, zhao2015bayesian}, brain computer interface \cite{mocks1988topographic}, signal processing \cite{de2008blind, muti2005multidimensional}, etc.

The main concept of solving tensor completion problem is that we use the observed entries of incomplete data to find the tensor decomposition factors which contain the latent features of the data, then we use the powerful feature representation ability of tensor decomposition factors to approximate the missing entries. The most studied and popular decomposition models in recent years are CANDECOMP/PARAFAC (CP) decomposition \cite{harshman1970foundations} and Tucker decomposition \cite{tucker1966some}. They have been applied in many data completion methods. CP weighted optimization (CP-WOPT) \cite{acar2011scalable} builds objective function by the Frobenius norm of weighted approximated tensor and observed tensor, then it uses optimization method to find the optimal CP factor matrices by the observed data.  Bayesian CP factorization \cite{zhao2015bayesian} employs Bayesian probabilistic model to find the best CP factor matrices and determine the rank of CP tensor automatically at the same time. The method in \cite{gandy2011tensor} recovers low-n-rank tensor data with its convex relaxation by alternating direction method of multipliers (ADM). Low-n-rank Tucker completion method is used in \cite{filipovic2015tucker} and the experiments show better results than other nuclear norm minimization methods.

Though CP and Tucker can reach relatively high performance in low-order tensors, due to the nature limitations of CP and Tucker, when it comes to high-order tensors and high missing rate of data, the performance of these two decomposition methods will decrease rapidly. Tensor-train (TT) decomposition \cite{oseledets2011tensor}, which is free from the ``curse of dimensionality'' and a better model to process high-order tensor is employed in our method. The works in our paper are concluded as below: (a) We propose an algorithm named Sparse Tensor-train Optimization (STTO) which considers incomplete data as sparse tensor and optimize the factors of tensor-train decomposition by gradient descent method. By optimizing the factors of tensor-train decomposition in sparse format, the computational complexity is significantly reduced. The tensor decomposition factors are used to approximate the missing entries. (b) Using synthetic data, we conduct simulation experiments to compare our algorithm with the state-of-the-art algorithms in four different dimensions. (c) We provide a data dimension ascending scheme for image data which can improve the performance of our algorithm. It is particularly useful to process image data in irregular missing cases like whole row missing and block missing. (d) We carry out several real world data experiments, and the results in simulation data and image data show that our method outperforms the state-of-the-art approaches. 

\section{Notations and Tensor-train Decomposition}

\subsection{Notations}
In this paper, we adopt the notations from \cite{kolda2009tensor}. Scalars are denoted by normal lowercase letters, e.g., $x$, and vectors are denoted by boldface lowercase letters, e.g., $\vect{x}$. Matrices are denoted by boldface capital letters, e.g., $\mat{X}$. Tensors of order $N\geq 3$ are denoted by boldface Euler script letters, e.g., $\tensor{X}$. $\mat{X}^{(n)}$ denotes the $n$th matrix of a matrix sequence, and the representations of vector and tensor sequence are denoted in the same way. When given a tensor $\tensor{X}  \in\mathbb{R}^{I_1\times I_2\times\cdots \times I_N}$, the $(i_{1},i_{2},\cdots,i_{N})$th element of $\tensor{X}$ is denoted by $x_{i_{1}i_{2}\cdots i_{N}}$ or $\tensor{X}(i_{1},i_{2},\cdots,i_{N})$.
 
 The $inner \ product$ of two tensors $\tensor{X}$, $\tensor{Y}\in\mathbb{R}^{I_1\times I_2\times\cdots \times I_N}$ is defined as $\langle \tensor{X},\tensor{Y} \rangle=\sum_{i_1}\sum_{i_2}\cdots\sum_{i_N}x_{i_1 i_2\cdots i_N}y_{i_1 i_2\cdots i_N}$. Furthermore, the $Frobenius \ norm$ of $\tensor{X}$ is defined by $\left \| \tensor{X} \right \|_F=\sqrt{\langle \tensor{X},\tensor{X} \rangle}$. The $Hadamard \ product$ is denoted by $\ast$ which is an element-wise product of vectors, matrices or tensors of same sizes. The $Kronecker \ product$ of two matrices $\mat{X}\in\mathbb{R}^{I \times K}$ and $\mat{Y}\in\mathbb{R}^{J \times L}$ is $\mat{X} \otimes \mat{Y} \in\mathbb{R}^{IJ \times KL}$. 
 
\subsection{Tensor-train Decomposition}
The most prominent advantage of tensor-train decomposition is that the amount of model parameters will not grow exponentially by data dimension. It decomposes a tensor into a sequence of three-way tensor factors (core tensors). In particular, the TT decomposition of a tensor $\tensor{X}  \in\mathbb{R}^{I_1\times I_2\times\cdots \times I_N}$ can be expressed as follow:
\begin{equation}
\label{tt_decom}
\tensor{X}=\ll \tensor{G}^{(1)},\tensor{G}^{(2)},\cdots,\tensor{G}^{(N)} \gg,
\end{equation}
where  $\tensor{G}^{(1)},\tensor{G}^{(2)},\cdots,\tensor{G}^{(N)} $ is a sequence of three-way core tensors of size $r_0 \times I_{1} \times r_{1},r_{1} \times I_{2} \times r_{2}, \cdots , r_{N-1} \times I_{N}\times r_N$, $r_0=r_N=1$. $\vect{r}= \{r_{0}, r_{1},r_{2},\cdots ,r_{N-1},r_{N}\}$ is named TT-rank which limits the size of every core tensor. Furthermore, Each element of tensor $\tensor{X}$ can be represented by core tensors as follow:
\begin{equation}
\label{TT_index}
 x_{i_{1}i_{2}\cdots i_{N}}=\prod\limits_{n=1}^{N}\mat{G}^{(n)}_{i_{n}},
\end{equation} 
where $\mat{G}^{(n)}_{i_{n}}$ is the $i_{n}$th slice of the $n$th core tensor of size $r_{n-1}\times r_{n}$, $n=1,2,\cdots,N$, $i_n\in\{ 1, 2, \cdots, I_n\}$. 

\section{Sparse Tensor-train Optimization}

\subsection{Our Previous Work}
In our previous work \cite{yuan2017completion}, we proposed an algorithm called Tensor-train Weighted OPTimization (TT-WOPT) which achieves high performance in data completion task. However, TT-WOPT considers all the missing entries of data as zero, and it computes the whole scale of tensor in every iteration. If the data scale is huge and missing rate is high, TT-WOPT will cost much computer memory space and be ineffective as it computes the whole scale tensor of which only a small percentage of entries is useful. 

\subsection{STTO Algorithm}
In order to solve the problems of TT-WOPT as mentioned in $Section\ 3.1$, our proposed algorithm STTO, which only uses observed entries to compute the gradient of every core tensor is proposed. Consider $\tensor{Y}$ is the observed tensor with missing entries, $\tensor{X}$ is the tensor approximated by core tensors, and the number of all the observed entries is $M$. Define the index of the $m$th observed entry as $\{ i_{1}^m, i_{2}^m,\cdots, i_{N}^m\}$, $m=1,\cdots,M$, we have $y_m=\tensor{Y}(i_{1}^m,i_{2}^m,\cdots,i_{N}^m)$, $x_m=\tensor{X}(i_{1}^m,i_{2}^m,\cdots,i_{N}^m)$. According to equation (\ref{TT_index}), $x_m$ can be written as:
\begin{equation}
\label{rec}
x_m=\prod\limits_{n=1}^{N}\mat{G}^{(n)}_{i_{n}^m}.
\end{equation}
For one observed entry of tensor $\tensor{Y}$, we formulate the objective function as: 
\begin{equation}
\label{so}
f(\mat{G}^{(1)}_{i_1^m}, \mat{G}^{(2)}_{i_2^m}, \cdots, \mat{G}^{(N)}_{i_N^m})=\frac{1}{2}\left\|y_m-\prod\limits_{n=1}^{N}\mat{G}^{(n)}_{i_{n}^m}\right\|^2_F.
\end{equation}
For $n=1,2,\cdots,N$, and $m=1,\cdots,M$, the partial derivatives of every used slice $\mat{G}^{(n)}_{i_n^m} $ of this entry is calculated by:
\begin{equation}
\label{ym}
\frac{\partial{f}}{\partial{\mat{G}^{(n)}_{i_n^m} }}=(x_m-y_m)(\mat{G}^{>n}_{i_n^m}\mat{G}^{<n}_{i_n^m})^T,
\end{equation}
where $\mat{G}^{>n}_{i_n^m}=\prod\limits_{n=n+1}^{N}\mat{G}^{(n)}_{i_{n}^m}$, $\mat{G}^{<n}_{i_n^m}=\prod\limits_{n=1}^{n-1}\mat{G}^{(n)}_{i_{n}^m}$.
If we consider the incomplete tensor as a sparse tensor, only the observed entries need to be enumerated. We arrange all the observed entries into vector $\vect{y}\in\mathbb{R}^M$, and arrange the according entries which are approximated by core tensors into $\vect{x}\in\mathbb{R}^M$. 
Then the optimization objective function of all missing entries can be formulated by:
\begin{equation}
f(\tensor{G}^{(1)},\tensor{G}^{(2)},\cdots,\tensor{G}^{(N)})=\frac{1}{2}\left\|\vect{y}-\vect{x}\right\|^2_F.
\end{equation}
By equation (\ref{rec}) and (\ref{so}), the optimization objective function can also be formulated as follow:
\begin{equation}
f(\tensor{G}^{(1)},\tensor{G}^{(2)},\cdots,\tensor{G}^{(N)})=\frac{1}{2}\sum^M_{m=1}\left\|y_m-x_m\right\|^2_F.
\end{equation}
So the sum gradient of every slice $\mat{G}^{(n)}_{j}$ of every core tensor is the accumulation of the slice gradients in equation (\ref{ym}) with the same index, that is:
\begin{equation}
\frac{\partial{f}}{\partial{\mat{G}^{(n)}_{j}}}=\sum^M_{{\substack{ m=1\\ m: i_n^m=j}}}(x_m-y_m)(\mat{G}^{>n}_{i_n^m}\mat{G}^{<n}_{i_n^m})^T,
\end{equation}
$j=1,2,\cdots,I_n$, and $n=1,2,\cdots,N$. After all the gradients of every slice of core tensors are obtained, any first-order optimization method can be applied to the STTO algorithm. The whole process of STTO is summarized in $Algorithm\ 1.$ The computational complexity of TT-WOPT and STTO is $\tensor{O}(r^{N-1}I^{N-1})$ and $\tensor{O}(Mr^{N-1})$, respectively. From this we can see STTO largely reduces the computational complexity and is totally free from dimensionality of tensor.

\begin{table}[H]
\footnotesize
\begin{center}
\begin{tabular}{l}
\hline
\textbf{Algorithm 1} Sparse Tensor-train Optimization (STTO)\\
\hline
1:\;\; \textbf{Input}: incomplete sparse tensor $\tensor{Y}$ and TT-rank $\vect{r}$.\\
2:\;\; \textbf{Initialization}: core tensors $\tensor{G}^{(1)},\tensor{G}^{(2)},\cdots,\tensor{G}^{(N)} $of approximated\\  \; \ \quad \quad  \quad \quad \quad  \quad \quad tensor $\tensor{X}$.\\
3:\;\; \textbf{While} the optimization stopping condition is not satisfied\\
4:\;\; \textbf{For} n=1:$N$\\
5:\;\;\;\; \textbf{For} j=1:$I_n$\\
6:\; \;\;\;\;\;Compute $\frac{\partial{f}}{\partial{\mat{G}^{(n)}_{j}}}=\sum^M_{{\substack{ m=1\\ m: i_n^m=j}}}(x_m-y_m)(\mat{G}^{>n}_{i_n^m}\mat{G}^{<n}_{i_n^m})^T$.\\
7:\;\;\;\; \textbf{End}\\
8:\;\; \textbf{End}\\
9:\;\; Update $\tensor{G}^{(1)},\tensor{G}^{(2)},\cdots,\tensor{G}^{(N)} $ by gradient descent method. \\
10: \textbf{End while} \\
11: \textbf{Output}: $\tensor{G}^{(1)},\tensor{G}^{(2)},\cdots,\tensor{G}^{(N)} $.\\
\hline
\end{tabular}
\end{center}
\end{table}

\section{Experiments}
In this section, our proposed STTO is compared with two state-of-the-art algorithms: CP weighted optimization (CP-WOPT) \cite{acar2011scalable} and Fully Bayesian CP (FBCP) \cite{zhao2015bayesian}. Simulation experiments, color image data experiments are conducted to validate the effectiveness of our algorithm. In addition, we provide a tensorization method to transform visual data to a higher dimension. This method can enhance the structure relation information of data and improve the performance of our algorithm. 

For evaluation indices, we use RSE (Relative Square Error) for simulation data and image data. PSNR (Peak Signal-to-noise Ratio) is used to measure the quality of reconstructed image data. In order to have a more clear comparison with CP-WOPT, we adopt the same optimization method as paper \cite{acar2011scalable}. We apply nonlinear conjugate gradient (NCG) with Hestenes-Stiefel updates \cite{nocedal2006numerical} and the Mor\'e-Thuente line search algorithm \cite{more1994line}. All the methods are implemented by an optimization toolbox named Pablano Toolbox \cite{dunlavy2010poblano}  and optimization stopping condition is set as maximum number of iterations.
\subsection{Simulations}
We consider to use values produced from a highly oscillating function: $f(x)=sin\frac{x}{4}cos(x^2)$ \cite{khoromskij2015tensor} as simulation data, which is expected to be well approximated by all the tensor completion algorithms. The four tested data structures are $26 \times 26 \times 26$ (3D), $7 \times7 \times7 \times7 \times7$ (5D), $4 \times4 \times4 \times4 \times4 \times4 \times 4$ (7D), $3 \times3 \times3 \times3 \times3 \times3 \times3 \times3 \times3$ (9D). The TT-ranks and CP-ranks of the four simulation are set to make the number of model parameters of the three algorithms as close as possible respectively.

From $Fig.\ 1.$ we can see, our method performs best among the three algorithms almost in every situation. Especially when the dimension of data is increase, our algorithm can maintain the RSE values while the performance of the other two algorithms falls quickly. 
\begin{figure}[ht]
\begin{center}
\label{sim}
\includegraphics[width=0.8\linewidth]{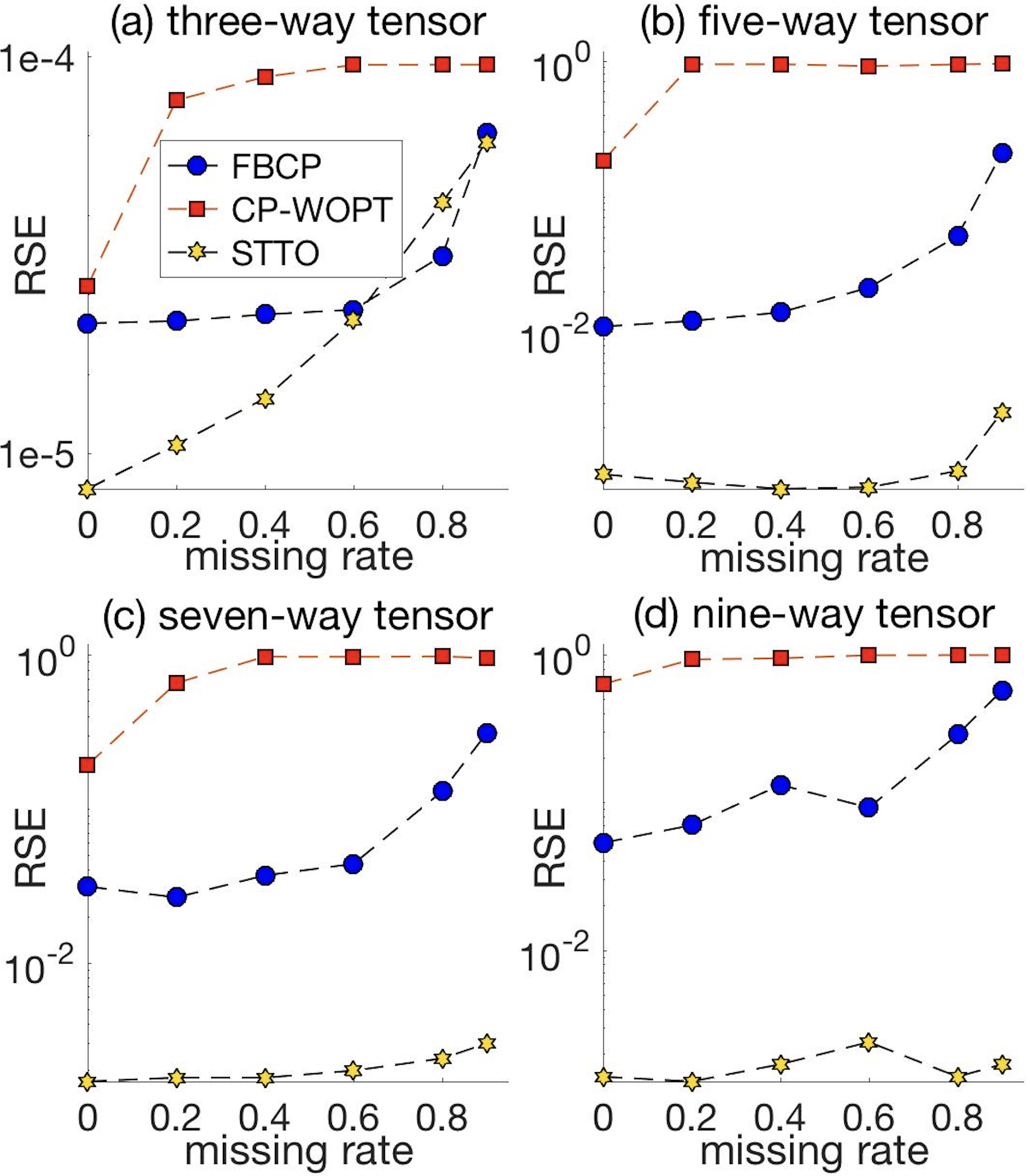}
\caption{RSE comparison of three algorithms under four different tensor dimensions. Missing rates of data are tested from 0\% to 90\%.}
\end{center}
\end{figure}

\subsection{Image Data Completion}
\subsubsection{Visual Data Tensorization Method}
From the simulation results we can see STTO can perform well in high-order cases, so we provide the below method to transform visual data to higher-order to enhance the performance of our algorithm. The original size of every image data is $256\times 256 \times 3$. First the three-way tensor image is reshaped to a seventeen-way tensor of size $2 \times2 \times \cdots \times2 \times3 $ and permute the tensor according to order $\{1 \;9 \;2 \;10 \;3 \;11 \;4 \;12 \;5 \;13 \;6 \;14 \;7 \;15 \;8 \;16 \;17\}$. Then we reshape the tensor to a nine-way tensor of size $4 \times4 \times\cdots  \times4 \times3$. The first order of the transformed tensor contains the data of a $2 \times 2$ pixel block of the image and the following orders of the tensor describe the expanding pixel blocks of the image. This nine-way tensor is considered to be a better structure of the image data. This tensorization method is applied to STTO in all of the following image experiments. The other two algorithms use original three way tensor form because they perform better in low-order tensor. 
\subsubsection{Random Missing}
We first adopt one benchmark image named ``Lena'' to see the best performance of all the algorithms in random missing cases. Briefly, we only compare the three algorithms in high missing rate situations. TT-ranks and CP-ranks are set properly to obtain the best results. The visualized experiment results in $Fig.\ 2.$ show that our STTO algorithm outperforms other algorithms distinctly. Particularly, when the missing rate reaches 98\% and 99\%, our algorithm with our visual data tensorization method can recover the image well while other algorithms fail totally.

\begin{figure}[htb]
\begin{center}
\includegraphics[width=1\linewidth]{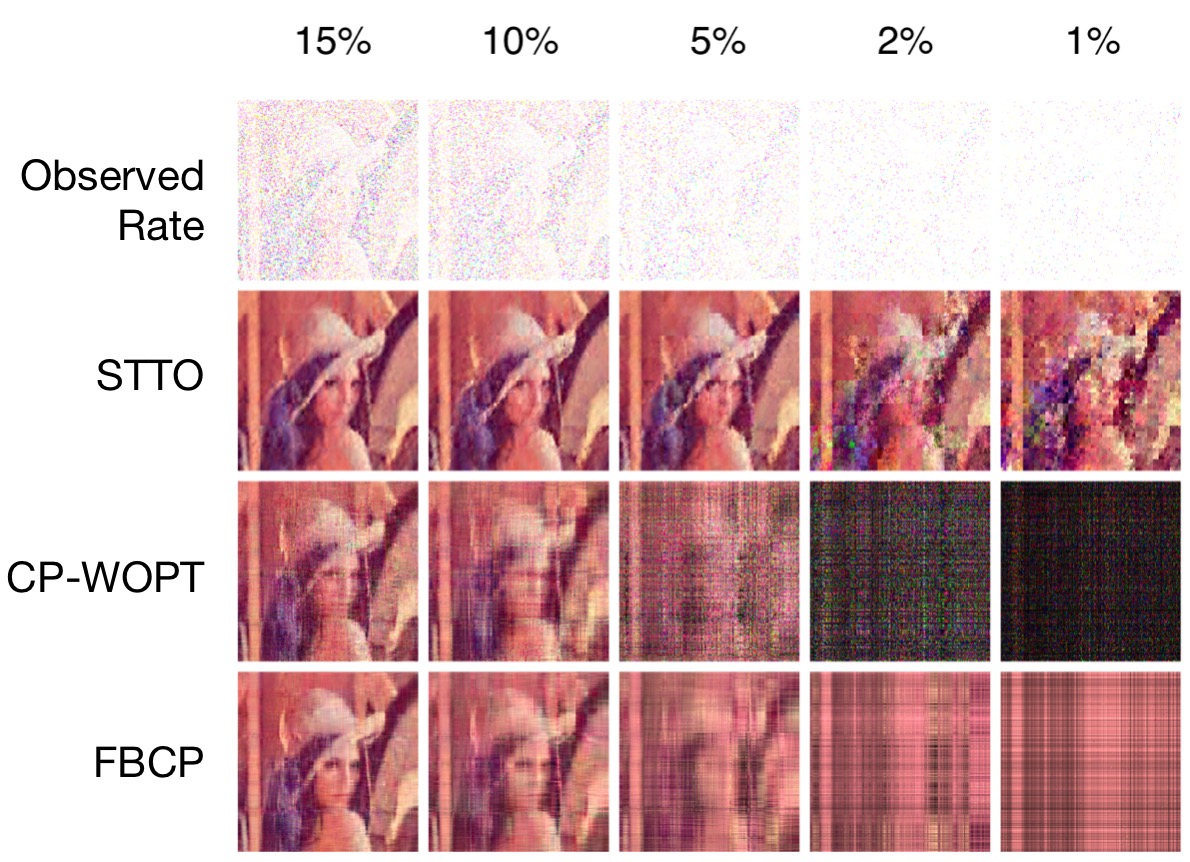}
\caption{Visualizing results of image recover performance of three algorithms under five missing rates.}
\label{rand}
\end{center}
\end{figure}
\subsubsection{Irregular Missing}
In this experiment, images with whole row missing or block missing are tested by the three algorithms. The visualized results of $Fig.\ 3.$ and values of RSE and PSNR from $Table\ 1.$ show that STTO with visual data tensorization method can recover images with whole row missing and block missing well.
\begin{figure}[h]
\begin{center}
\includegraphics[width=1\linewidth]{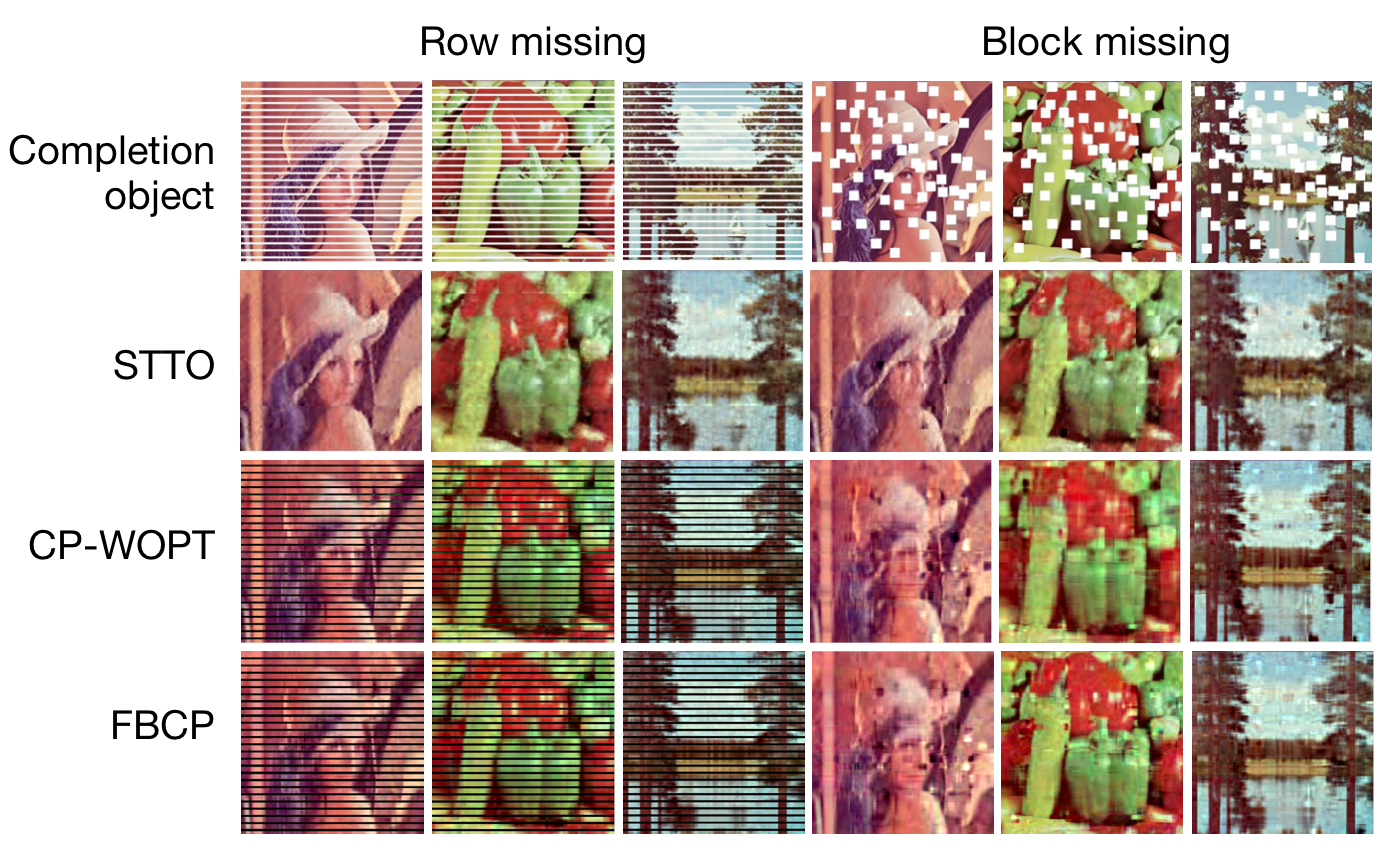}
\caption{Visualizing results of image recover performance of three different algorithms under two special missing cases.}
\label{irregu}
\end{center}
\end{figure}
\setlength\tabcolsep{1.5pt} 
\begin{table}[h]\footnotesize
\centering
\caption{Comparison of the recover performance (RSE and PSNR) of three algorithms under two special missing cases. }
\label{testing_results}
\begin{tabular}{c|c|c|c|c|c|c|c}
\hline
\hline
 \multicolumn{2}{c|}{}& \multicolumn{3}{c|}{row missing}& \multicolumn{3}{c}{block missing}\\
 \hline
 \multicolumn{2}{c|}{image}& lena & peppers & sailboat & lena & peppers&sailboat\\
 \hline
STTO&\makecell[cc]{RSE \\ PSNR}   &     \makecell[cc]{\textbf{0.1138} \\ \textbf{24.00}}    &    \makecell[cc]{\textbf{0.1661} \\ \textbf{20.80}}     &
\makecell[cc]{\textbf{0.1767}\\ \textbf{19.93}}     &           \makecell[cc]{\textbf{0.1323}\\ \textbf{22.69}}   &\makecell[cc]{\textbf{0.1611}\\ \textbf{21.06}}&\makecell[cc]{\textbf{0.1704}\\ \textbf{20.25}}\\
\hline
CP-WOPT&\makecell[cc]{RSE \\ PSNR}&\makecell[cc]{0.5401\\ 10.86}&\makecell[cc]{0.5546 \\ 10.85}&\makecell[cc]{0.5545 \\ 10.34}&\makecell[cc]{0.1746 \\20.61} &   \makecell[cc]{0.2252\\18.27}&\makecell[cc]{0.2082\\ 19.00}\\
\hline
FBCP&\makecell[cc]{RSE \\ PSNR}&\makecell[cc]{0.5503\\ 10.46}&\makecell[cc]{0.5594 \\ 10.58}&\makecell[cc]{0.5586\\ 10.18}&\makecell[cc]{0.1498 \\ 21.66}  & \makecell[cc]{0.1671\\ 20.79}&\makecell[cc]{0.1764\\ 20.01}\\
\hline\hline
\end{tabular}
\end{table}

\section{Conclusions}
In this paper, we first elaborate the basis of tensor and tensor-train decomposition. Then STTO algorithm which is efficient and has low computational complexity is proposed. It uses observed entries of sparse tensor to optimize the core tensors of tensor-train model and recover the missing data. From the simulation experiments, we can see our algorithm outperforms the state-of-the-art methods in both low-order cases and high-order cases. In addition, image completion experiment results prove that STTO with our tensorization scheme can achieve a high performance under high missing rate cases. The remarkable results on image irregular missing cases also show advantages of our algorithm and tensorization method. From the experiment results we can see tensor-train decomposition with high-order tensorizations can achieve high compression and representation abilities. Furthermore, it should be noted that the performance of tensor-train decomposition is sensitive to the selection of TT-ranks. Hence, we will study on how to optimize tensor factors and TT-ranks simultaneously in our future work. 

\section{Acknowledgement}
This work was supported by JSPS KAKENHI (Grant No. 17K00326, 15H04002), JST CREST (Grant No. JPMJCR1784) and the National Natural Science Foundation of China (Grant No. 61773129).

\vfill\pagebreak 

\bibliographystyle{IEEEbis}
\bibliography{paper}

\end{document}